\documentclass[10pt]{article}
\usepackage{stmaryrd}
\usepackage{amsfonts}
\usepackage{bbm}

\usepackage{mathrsfs}
\usepackage{latexsym,amssymb,amsmath,amscd,amsthm,amsxtra,xypic}
\usepackage[dvips]{graphicx}

\usepackage{color}
\usepackage{xcolor}
\definecolor{red}{rgb}{1,0,0}

\newtheorem{thm}{Theorem}[section]
\newtheorem{lem}[thm]{Lemma}
\newtheorem{cor}[thm]{Corollary}
\newtheorem{pro}[thm]{Proposition}
\newtheorem{ex}[thm]{Example}
\newtheorem{rmk}[thm]{Remark}
\newtheorem{defi}[thm]{Definition}

\newcommand {\comment}[1]{{\marginpar{*}\scriptsize\textbf{Comments:} #1}}

\newcommand{\Dorfman}[1]{\left  \llbracket  #1\right  \rrbracket }

\newcommand{\CDorfman}[1]{   [    #1   ]_{\huaK}   }
\newcommand{\pair}[1]{\left( #1\right)_+}
\newcommand{\pairE}[1]{\left ( #1\right )_E}
\newcommand{\pairV}[1]{\left ( #1\right )_V}
\newcommand{\ppE}[1]{\left\langle  #1\right\rangle }
\newcommand{\pe}[1]{\left \langle #1\right\rangle_E}
\newcommand{\pibracket}[1]{\left [ #1\right ]_{\pi}}
\newcommand{\pairL}[1]{\left( #1\right)_{L^\circ}}

\newcommand{\LN }{\mathsf{LN}}
\newcommand{\JN }{\mathsf{JN}}
\newcommand{\VB }{\mathsf{VB}}
\newcommand{\E }{\mathsf{E}}
\newcommand{\pr }{\mathsf{pr}}
\newcommand{\ol }{\frko\frkl}
\newcommand{\g }{\frkg}
\newcommand {\yh}[1]{{\marginpar{*}\scriptsize\textcolor{purple}{yh: #1}}}
\newcommand{\Id}{\mathsf{id}}

\newcommand {\emptycomment}[1]{}

\newcommand{\lon }{\,\rightarrow\,}
\newcommand{\be }{\begin{eqnarray*}}
\newcommand{\ee }{\end{eqnarray*}}

\newcommand{\defbe}{\triangleq}
\newcommand{\pf}{\noindent{\bf Proof.}\ }

\newcommand{\huaA}{\mathcal{A}}
\newcommand{\huaL}{\mathcal{L}}

\newcommand{\huaE}{\mathcal{E}}
\newcommand{\huaF}{\mathcal{F}}

\newcommand{\huaC}{{\mathcal{C}}}

\newcommand{\huaJ}{\mathcal{J}}

\newcommand{\huaK}{\mathcal{K}}

\newcommand{\set}[1]{\left\{#1\right\}}

\newcommand{\frkd}{\mathfrak d}

\newcommand{\frkg}{\mathfrak g}

\newcommand{\frkk}{\mathfrak k}
\newcommand{\frkl}{\mathfrak l}
\newcommand{\frkm}{\mathfrak m}

\newcommand{\frko}{\mathfrak o}

\newcommand{\frkr}{\mathfrak r}

\newcommand{\frkt}{\mathfrak t}

\newcommand{\frky}{\mathfrak y}

\newcommand{\frkL}{\mathfrak L}

\newcommand{\frkX}{\mathfrak X}

\def\gpd{\,\lower1pt\hbox{$\longrightarrow$}\hskip-.24in\raise2pt
         \hbox{$\longrightarrow$}\,}

\def\qed{\hfill ~\vrule height6pt width6pt depth0pt}

\newcommand{\Lied}{\frkL}

\newcommand{\pisharp}{\pi}

\newcommand{\half}{\frac{1}{2}}

\newcommand{\jet}{\mathfrak{J}}

\newcommand{\jetd}{\mathbbm{d}}
\newcommand{\dev}{\mathfrak{D}}

\newcommand{\dA}{\mathrm{d}^A}
\newcommand{\dB}{\mathrm{d}^B}

\newcommand{\pie}{^\prime}

\newcommand{\e}{\mathbbm{e}}
\newcommand{\p}{\mathbbm{p}}
\newcommand{\id}{\mathbbm{i}}
\newcommand{\jd}{\mathbbm{j}}

\newcommand{\dM}{\mathrm{d}}

\newcommand{\Hom}{\mathrm{Hom}}
\newcommand{\gl}{\mathrm{gl}}

\newcommand{\rhowx}{\rho^{\star}}

\newcommand{\Img}{\mathrm{Im}}
\newcommand{\ad}{\mathrm{ad}}

\begin{document}
\title{
{ $\VB$-Courant algebroids, $\E$-Courant algebroids and generalized geometry
\thanks
 {
 Research supported by NSFC (11471139) and NSF of Jilin Province (20170101050JC).
 }
} }
\author{\vspace{1mm}Honglei Lang$^1$, Yunhe Sheng$^2$ and Aissa Wade$^3$ \\
$^1$Max Planck Institute for Mathematics,\\\vspace{1mm}
Bonn 53111, Germany\\
$^2$Department of Mathematics, Jilin University,\\
 Changchun 130012, Jilin, China\vspace{1mm}
 \\
 $^3$Mathematics Department, Penn State University,\\\vspace{1mm} University Park, PA 16802, USA. \\
email:  hllang@mpim-bonn.mpg.de,~~
shengyh@jlu.edu.cn,~~ wade@math.psu.edu
 }

\date{}


\maketitle

\begin{abstract}
In this paper, we  first discuss the relation between $\VB$-Courant algebroids and $\E$-Courant algebroids and construct some examples of $\E$-Courant algebroids. Then we introduce the notion of  a generalized complex
structure on an $\E$-Courant algebroid, unifying  the usual
generalized complex structures on  even-dimensional manifolds and
generalized contact structures on  odd-dimensional manifolds. Moreover, we study generalized complex structures on an omni-Lie algebroid in detail. In particular, we show that generalized complex structures on an omni-Lie algebra $\gl(V)\oplus V$ correspond to complex Lie algebra structures on $V$.

\end{abstract}


\section{Introduction}
The theory of Courant algebroids was first  introduced by  Liu, Weinstein and Xu \cite{LWXmani} providing an extension of Drinfeld's double for Lie bialgebroids. Precisely, the double of a Lie bialgebroid is a special Courant algebroid \cite{LWXmani,Roytenbergphdthesis}. Jacobi algebroids are natural extensions of Lie algebroids  and Courant-Jacobi algebroids. They were considered by Grabowski and Marmo \cite{Grabowski marmo2}  and they can be viewed as generalizations of Courant algebroids.
Both Courant algebroids and Courant-Jacobi algebroids have been extensively studied in the last decade since these are crucial geometric tools in  Poisson geometry and mathematical physics. It is known that they both belong to a more general framework, namely that of  $\E$-Courant algebroids.  Indeed,  $\E$-Courant algebroids were introduced by Chen, Liu and the second author in \cite{CLS2} as a differential geometric  object encompassing  Courant algebroids \cite{LWXmani}, Courant-Jacobi algebroids \cite{Grabowski marmo2}, omni-Lie algebroids \cite{CLomni}, conformal Courant algebroids \cite{conformalCA} and $AV$-Courant algebroids \cite{AV courant}. It turns out that $\E$-Courant algebroids are related to more geometric structures such as $\VB$-Courant algebroids \cite{Li-Bland}.
\medskip

The aim of this  paper is two-fold. Firstly, we illuminate  the relationship between $\VB$-Courant algebroids and $\E$-Courant algebroids.
Secondly,  we study generalized complex structures on $\E$-Courant algebroids.  Recall that a generalized almost complex structure on a manifold $M$ is an endomorphism ${\mathcal J}$ of the \emph{generalized tangent bundle $\mathbb T M := TM \oplus T^\ast M$}  which preserves the natural pairing on $\mathbb T M$ and such that ${\mathcal J}^2 = - \operatorname{id}$. If, additionally,  the $\sqrt{-1}$-eigenbundle of ${\mathcal J}$ in the complexification $\mathbb TM \otimes \mathbb C$  is involutive relative to the Dorfman (equivalently, the Courant) bracket, then ${\mathcal J}$ is said to be \emph{integrable}, and  $(M, \mathcal J)$ is called a \emph{generalized complex manifold}. See \cite{BCG,CrainicGCS,gualtieri,hitchin,SXreduction} for more details.

\medskip
Given  a vector bundle $E\stackrel{q}{\lon}M$, we consider  its gauge Lie algebroid $\dev E$, i.e. the gauge Lie algebroid of the
 frame bundle $\huaF(E)$. It is known that $\dev{E}$ is a  transitive Lie algebroid over $M$  and the first jet bundle $\jet E$ is its  $E$-dual bundle.  In fact,  $\ol(E)=\dev E\oplus \jet E$ is called an omni-Lie algebroid \cite{CLomni}, which is a generalization of  Weinstein's  concept of an  omni-Lie algebra \cite{weinstein:omni}. In particular, the  line bundle case where $E$ comes from a contact distribution brings us to the concept of a generalized contact bundle. To have a better grasp of the concept of a generalized contact bundle, we briefly review the line bundle approach to contact geometry. By definition,  a contact structure on an odd-dimensional manifold $M$ is a maximal non-integrable hyperplane distribution $H \subset TM$. In a dual way, any hyperplane distribution $H$ on $M$ can be regarded as a nowhere vanishing $1$-form $\theta : TM \to L$ {(its \emph{structure form})} with values in the line bundle $L =TM/ H$, {such that $H = \ker \theta$}. Replacing the tangent algebroid with the Atiyah algebroid of a line bundle in the definition of a generalized complex manifold, we obtain the notion of a \emph{generalized contact bundle}. In  this paper,  we extend the concept of a generalized contact bundle to the context of $\E$-Courant algebroids.

 \medskip

 The paper is organized as follows. Section 2  contains basic definitions used in the sequel. Section 3
highlights the importance and naturality of the notion of $\E$-Courant algebroids. Explicitly, {\it  the fat Courant algebroid associated to a    $\VB$-Courant algebroid }  (see  the definition of a $\VB$-Courant algebroid below) is an $\E$-Courant algebroid. We observe the following facts:

\begin{itemize}
\item[$\bullet$] Given a crossed module of Lie algebras $(\mathfrak{m},\mathfrak{g})$, we get an $\frkm$-Courant algebroid $\Hom(\mathfrak{g},\mathfrak{m})\oplus \mathfrak{g}$, which was given in \cite{LSX} as a generalization of an omni-Lie algebra.

\item[$\bullet$] The omni-Lie algebroid $\ol(E)=\dev E\oplus \jet E$    is the linearization of the $\VB$-Courant algebroid $TE^*\oplus T^*E^*$. This generalizes the fact that an omni-Lie algebra is the linearization of the standard Courant algebroid.

\item[$\bullet$] For a Courant algebroid $\huaC$, $T \huaC$ is a $\VB$-Courant algebroid. The associated fat Courant algebroid  $\jet\huaC$ is  a $T^*M$-Courant algebroid. The fact that  $\jet \huaC$ is a $T^*M$-Courant algebroid     was first obtained in \cite[Theorem 2.13]{CLS2}.

\end{itemize}
In Section 4, we introduce generalized complex structures on $\E$-Courant algebroids and provide examples.  In Sections 5, we describe generalized complex structures on omni-Lie algebroids. In Section 6, we show that generalized complex structures on the omni-Lie algebra $\ol(V)$ are in one-to-one correspondence  with  complex Lie algebra structures on $V$.

\section{Preliminaries}
Throughout the paper, $M$ is a smooth manifold, $\dM$ is the usual
differential operator on forms, and $E\longrightarrow M$ is  a vector
bundle.  In this section, we recall the notions of $\E$-Courant
algebroids introduced in \cite{CLS2}, omni-Lie algebroids introduced in \cite{CLomni}, generalized
complex structures introduced   in \cite{gualtieri,hitchin} and
generalized contact structures introduced  in \cite{VitaglianoWade}.

\subsection{ $\E$-Courant algebroids and omni-Lie algebroids}

For a vector bundle $E{\lon}M$, its gauge Lie algebroid $\dev E$
with the commutator bracket $[\cdot,\cdot]_\dev$ is just the gauge Lie algebroid
of the
 frame bundle
 $\huaF(E)$, which is also called the covariant differential operator bundle of $E$ (see \cite[Example
3.3.4]{Mkz:GTGA}). The corresponding Atiyah sequence is as
follows:
\begin{equation}\label{Seq:DE}
\xymatrix@C=0.5cm{0 \ar[r] & \gl(E)  \ar[rr]^{\id} &&
                \dev{E}  \ar[rr]^{\jd} && TM \ar[r]  & 0.
                }
\end{equation}
In \cite{CLomni},
the authors proved that the jet bundle $\jet E$   may be considered as
an $E$-dual bundle of $\dev E$:
\begin{equation}\label{eqn:jetE}  {\jet E} \cong
\set{\nu\in \Hom(\dev{E},E)\,|\,
\nu(\Phi)=\Phi\circ\nu(\Id_{E}),\quad\forall~ ~ ~\Phi\in \gl(E)}.
\end{equation}
Associated to the jet bundle $\jet E$,  there is a jet sequence
given by:
\begin{equation}\label{Seq:JetE}
\xymatrix@C=0.5cm{0 \ar[r] & \Hom(TM,E)
\ar[rr]^{\quad\quad\quad\e} &&
                {\jet}{E} \ar[rr]^{\p} && E \ar[r]  & 0.
                }
\end{equation}
Define the operator $\jetd: \Gamma(E) \rightarrow \Gamma(\jet E)$
 by:
$$\jetd u (\frkd) := \frkd (u),\quad \quad \forall~  u  \in   \Gamma(E),  ~~~\frkd \in\Gamma(\dev E).$$
An important formula which will be often used is
\begin{equation*}
\jetd(fu)=\dM f\otimes u+f\jetd
u,\quad \forall~u\in\Gamma(E),~f\in C^\infty(M).
\end{equation*}
In fact, there is an $E$-valued pairing between $\jet{E}$ and $\dev{E}$
by setting:
\begin{equation}\label{Eqt:conpairingE}
\ppE{\mu,\frkd}_E~ \defbe
\frkd(u),\quad\forall~ ~~ \mu\in (\jet{E})_m,~\frkd\in(\dev{E})_m,
\end{equation}
where $u\in \Gamma(E)$ satisfies $\mu=[u]_m$. In particular, one has
\begin{eqnarray*}
\ppE{\mu,\Phi}_E &=& \Phi\circ \p(\mu),\quad\forall~ ~~ \Phi\in \gl(E),~\mu\in\jet{E};\\
\label{conpairing2} \ppE{{\frky},\frkd}_{E} &=& {\frky}\circ
\jd(\frkd),\quad\forall~ ~~ \frky\in \Hom(TM,E),~\frkd\in\dev{E}.
\end{eqnarray*}

For vector bundles $P,Q$ over $M$ and a bundle map
$\rho:P\longrightarrow Q$, we denote the induced
$E$-dual bundle map by $\rho^\star$, i.e.
$$
\rho^\star : \Hom(Q,E)\rightarrow \Hom(P,E), \quad
\rho^\star(\nu)(k)=\nu(\rho(k)), \quad \forall~k\in
P,~\nu\in\Hom(Q,E).
$$

\begin{defi}{\rm(\cite{CLS2})}\label{defi: $E$-C A}
An \emph{ $\E$-Courant algebroid} is a quadruple $(\huaK,\CDorfman{\cdot,\cdot},
\pairE{\cdot,\cdot},\rho)$, where
$\huaK$ is a vector bundle  over $M$ such that
$(\Gamma(\huaK),\CDorfman{\cdot,\cdot})$ is a Leibniz algebra,
$\pairE{\cdot,\cdot}: ~\huaK\otimes\huaK\lon E$ is a
nondegenerate symmetric $E$-valued pairing, which induces an
embedding: $\huaK \hookrightarrow \Hom(\huaK,E)$ via $Y(X)=2\pairE{X,Y}$, and
$\rho:~ \huaK\lon \dev E$ is a bundle map called the anchor, such that for all $~X,Y,Z\in \Gamma(\huaK)$, the following
properties hold:
\begin{itemize}
\item[~~~~ \rm (EC-1)]  \quad \quad \quad \quad
$\rho\CDorfman{X,Y} =[\rho(X),\rho(Y)]_{\dev}$;
\item[~~~~ \rm (EC-2)] \quad \quad \quad \quad  $\CDorfman{X,X}  =  \rhowx\jetd\pairE{X,X}$;
\item[~~~~ \rm (EC-3)]\quad \quad \quad \quad
$\rho(X)\pairE{Y,Z}=\pairE{\CDorfman{X,Y} ,Z}
+\pairE{Y,\CDorfman{X,Z} } $;
 \item[~~~~ \rm (EC-4)]  \quad \quad \quad \quad $\rhowx( \jet E)\subset \huaK$,
i.e. $ \pairE{\rhowx(\mu), X}=\half\mu(\rho(X)),~ ~~\forall~ \mu\in\jet E;$
\item[~~~~ \rm (EC-5)] \quad \quad \quad \quad
$\rho\circ\rhowx=0$.
\end{itemize}
\end{defi}

Obviously, a Courant algebroid is an $\E$-Courant algebroid, where $E=M\times \mathbb R$, the trivial line bundle. Similar to the proof for Courant algebroids (\cite[Lemma 2.6.2]{Roytenbergphdthesis}), we have
\begin{lem}\label{[e,Df]}
For an $\E$-Courant algebroid $\huaK$, one has
\[[X,\rhowx\jetd u]_{\huaK}=2\rhowx\jetd \pairE{X,\rhowx\jetd u},\quad \quad [\rhowx\jetd u,X]_{\huaK}=0,\quad \quad \forall~X\in \Gamma(\huaK), ~u\in \Gamma(E).\]
\end{lem}

An omni-Lie algebroid, which was introduced in
\cite{CLomni}, is a very interesting example of $\E$-Courant algebroids. Let us recall it briefly.
There is an $E$-valued pairing
$\pairE{\cdot,\cdot}$ on $\dev E\oplus\jet E$ defined by
\begin{equation}\label{eqn:pair omni}
\pairE{\frkd+\mu,\frkt+\nu}=\frac{1}{2}\big(\ppE{\mu,\frkt}_E+\ppE{\nu,\frkd}_E\big),\quad\forall~\frkd+\mu,~\frkt+\nu\in\dev
E\oplus\jet E.
\end{equation}
Furthermore,  $\Gamma(\jet E)$ is invariant under the Lie derivative $\Lied_{\frkd}$ for any
 $\frkd \in\Gamma(\dev E)$
 which is defined by the
 Leibniz rule:
\begin{eqnarray}\nonumber
\ppE{\Lied_{\frkd}\mu,\frkd\pie}_{E}&\defbe&
\frkd\ppE{\mu,\frkd\pie}_{E}-\ppE{\mu,[\frkd,\frkd\pie]_{\dev}}_{E},
\quad\forall~ \mu \in \Gamma(\jet{E}), ~
~\frkd\pie\in\Gamma(\dev{E}).
\end{eqnarray}
On the section space $\Gamma(\dev E\oplus \jet E)$, we can define
a bracket as follows
\begin{eqnarray}\label{eqn:backet omni}
\Dorfman{\frkd+\mu,\frkr+\nu}&\defbe&
[\frkd,\frkr]_{\dev}+\Lied_{\frkd}\nu-\Lied_{\frkr}\mu +
\jetd\mu(\frkr).
\end{eqnarray}

\begin{defi}{\rm(\cite{CLomni})}
The quadruple
$(\dev E\oplus \jet
E,\Dorfman{\cdot,\cdot},\pairE{\cdot,\cdot},\rho)$ is called
an \emph{ omni-Lie algebroid},  where $\rho$ is the projection from
$\dev E\oplus \jet E$ to $\dev{E}$, $\pairE{\cdot,\cdot}$ and
$\Dorfman{\cdot,\cdot}$ are given by $(\ref{eqn:pair omni})$ and
$(\ref{eqn:backet omni})$ respectively.
\end{defi}

We will denote an omni-Lie algebroid by $\ol(E)$.

\subsection{ Generalized complex structures and generalized contact
structures}

The notion of a Courant algebroid was introduced in \cite{LWXmani}.  A Courant algebroid is a quadruple $(\huaC,\Dorfman{\cdot,\cdot},\pair{\cdot,\cdot},\rho)$, where $\huaC$ is a vector bundle over $M$, $\Dorfman{\cdot,\cdot}$ is a bracket operation on $\Gamma(\huaC)$, $\pair{\cdot,\cdot}$ is a nondegenerate symmetric bilinear form on $\huaC$, $\rho:\huaC\longrightarrow TM$ is a bundle map called the anchor, such that some compatibility conditions are satisfied. See \cite{Roytenbergphdthesis} for more details.
Consider the generalized tangent bundle
$$\mathbb T M:=TM\oplus T^*M.$$ On its
section space $\Gamma(\mathbb T M)$, there is a Dorfman
bracket:
\begin{equation}\label{eqn:bracket}
\Dorfman{X+\xi,Y+\eta}=[X,Y]+\huaL_X\eta-i_Y\dM\xi,\quad\forall~X+\xi,~Y+\eta\in\Gamma(\mathbb T M).
\end{equation}
 Furthermore, there is a canonical nondegenerate symmetric
bilinear form on $\mathbb T M$:
\begin{equation}\label{eqn:pair}
\pair{ X+\xi,Y+\eta}=\frac{1}{2}\big(\eta(X)+\xi(Y)\big).
\end{equation}
We call $(\mathbb T M,\Dorfman{\cdot,\cdot},\pair{\cdot,\cdot},\pr_{TM})$ the standard Courant algebroid.
\begin{defi}
A \emph{generalized complex structure} on a manifold $M$ is a bundle map
$\huaJ:\mathbb T M\longrightarrow\mathbb T M$ satisfying the algebraic
properties:
$$\huaJ^2=-\Id\quad\mbox{and}\quad\pair{\huaJ(u),\huaJ(v)}=\pair{ u,v}$$
and the integrability condition:
$$\Dorfman{\huaJ(u),\huaJ(v)}-\Dorfman{u,v}-\huaJ\big(\Dorfman{\huaJ(u),v}+\Dorfman{u,\huaJ(v)}\big)=0,\quad\forall~u,~v\in\Gamma(\mathbb T M).$$
Here, $\pair{\cdot,\cdot}$ and $\Dorfman{\cdot,\cdot}$ are given by
$(\ref{eqn:pair})$ and $(\ref{eqn:bracket})$ respectively.
\end{defi}

See \cite{gualtieri,hitchin} for more details. Note that only even-dimensional manifolds can have generalized
complex structures. In \cite{VitaglianoWade}, the
authors give the odd-dimensional analogue of the concept of a
generalized complex structures extending the definition given in  \cite{Wade contact manifold and g}.
We now recall the definition of a generalized contact bundle from \cite{VitaglianoWade}.
A \emph{generalized contact bundle} is a line bundle $L \to M$ equipped with a \emph{generalized contact structure}, i.e.~a vector bundle endomorphism
${\mathcal I}: \dev L \oplus\jet L \to \dev L\oplus\jet L$ such that
\begin{itemize}
\item ${\mathcal I}$ is \emph{almost complex}, i.e.~${\mathcal I}^2 = -\Id$,
\item ${\mathcal I}$ is \emph{skew-symmetric}, i.e.
\[
( {\mathcal I} \alpha, \beta)_L + ( \alpha, {\mathcal I} \beta )_L = 0, \quad \forall \alpha, \beta \in \Gamma ({\dev} L\oplus \jet L),
\]
\item ${\mathcal I}$ is \emph{integrable}, i.e.
\[
[\![ {\mathcal I} \alpha, {\mathcal I} \beta]\!] - [\![ \alpha, \beta]\!] - {\mathcal I} [\![ {\mathcal I} \alpha, \beta]\!] - {\mathcal I} [\![ \alpha, {\mathcal I} \beta ]\!] = 0, \quad \forall \alpha, \beta \in \Gamma ({\dev} L\oplus \jet L).
\]
\end{itemize}
Let $(L \to M, {\mathcal I})$ be a generalized contact bundle. Using the direct sum  $\ol(L)=\dev L \oplus\jet L$ and  the definition, one can see that
\[
{\mathcal I} = \left(
\begin{array}{cc}
\phi & J^\sharp \\
\omega_\flat & - \phi^\dag
\end{array}
\right),
\]
where $J$ is a Jacobi bi-derivation, $\phi$ is an endomorphism  of $\dev L$ compatible with $J$, and the $2$-form $\omega : \wedge^2 \dev L \to L$ and its associated vector bundle morphism $\omega_\flat : \dev L \to \mathfrak J L$ satisfy additional compatibility conditions \cite{VitaglianoWade}.

\section{$\VB$-Courant algebroids and $\E$-Courant algebroids}
In this section, we highlight the relation between $\VB$-Courant algebroids and $\E$-Courant algebroids and give more examples of $\E$-Courant algebroids.

Denote a double vector bundle
\begin{equation*}
\xymatrix{ D \ar@{->} [r]^{\pi_B}\ar@ {->} [d]_{\pi_A}
&B \ar@ {->} [d]^{q_B}\\
{A} \ar@{->} [r]_{q_A}&M }
\end{equation*}
with core $C$ by $(D;A,B;M)$.
The space of sections $\Gamma_B(D)$ is generated as a $C^\infty(B)$-module by core sections $\Gamma_B^c(D)$ and linear sections $\Gamma_B^l(D)$. See \cite{ME} for more details.
For a section $c: M\rightarrow C$, the corresponding {\bf core section} $c^\dag: B\rightarrow D$ is defined as
$$c^\dag(b_m)=\tilde{0}_{b_m}+_A \overline{c(m)},\quad \forall~ m\in M, b_m\in B_m,$$ where $\bar{\cdot}$ means the inclusion $C\hookrightarrow D$.
A section $\xi: B\rightarrow D$ is called \emph{linear} if it is a bundle morphism from $B\rightarrow M$ to $D\rightarrow A$ over a section $a\in \Gamma(A)$.
Given $\psi\in \Gamma(B^*\otimes C)$, there is a linear section $\tilde{\psi}: B\rightarrow D$ over the zero section $0^A: M\rightarrow A$ given by
$$\tilde{\psi}(b_m)=\tilde{0}_{b_m}+_A \overline{\psi(b_m)}.$$

Note that $\Gamma_B^l(D)$ is locally free as a $C^\infty (M)$-module. Therefore, $\Gamma_B^l(D)$ is equal to $\Gamma(\hat{A})$ for some vector bundle $\hat{A}\rightarrow M$. Moreover, we have the following short exact sequence of vector bundles over $M$
\begin{equation}\label{fat}
0 \rightarrow B^*\otimes C \longrightarrow \hat{A}\longrightarrow A \rightarrow 0.
\end{equation}

 \begin{ex}\label{fat seq}{\rm
 Let $E$ be a vector bundle  over $M$.
\begin{itemize}
  \item[\rm(a)]The tangent bundle $(TE;TM,E;M)$ is a double vector bundle with core $E$. Then $\hat{A}$ is   the gauge bundle $\dev E$ and the exact sequence (\ref{fat}) is exactly the Atiyah sequence (\ref{Seq:DE}).

  \item[\rm(b)]The cotangent bundle $(T^*E;E^*,E;M)$ is a double vector bundle with core $T^*M$. In this case, $\hat{A}$ is exactly the jet bundle $\jet E^*$  and the exact sequence (\ref{fat}) is indeed the jet sequence (\ref{Seq:JetE}).
\end{itemize}
}
\end{ex}
 \begin{defi} {\rm(\cite{Li-Bland})}
A \emph{$\VB$-Courant algebroid} is a metric double vector bundle
\begin{equation*}
\xymatrix{ \mathbb{E} \ar@{->} [r]^{\pi_B}\ar@ {->} [d]_{\pi_A}
&B \ar@ {->} [d]^{q_B}\\
{A} \ar@{->} [r]_{q_A}&M, }
\end{equation*}
with core $C$ such that $\mathbb{E}\rightarrow B$ is a Courant algebroid and the following conditions are satisfied:
\begin{itemize}
\item[\rm(i)] The anchor map $\Theta:\mathbb{E}\rightarrow TB$ is linear. That is
$$\Theta: (\mathbb{E};A,B;M)\longrightarrow (TB;TM,B;M)$$
is a morphism of double vector bundles.
\item[\rm(ii)] The Courant bracket is linear. That is
$$\Dorfman{\Gamma_B^l(\mathbb{E}),\Gamma_B^l(\mathbb{E})}\subseteq \Gamma_B^l(\mathbb{E}),\ \ \ \Dorfman{\Gamma_B^l(\mathbb{E}),\Gamma_B^c(\mathbb{E})}\subseteq \Gamma_B^c(\mathbb{E}),\quad \Dorfman{\Gamma_B^c(\mathbb{E}),\Gamma_B^c(\mathbb{E})}=0.$$
\end{itemize}
\end{defi}
For a $\VB$-Courant algebroid $\mathbb{E}$, we have the exact sequence $(\ref{fat})$. Note that the restriction of the pairing on $\mathbb{E}$ to linear sections of $\mathbb{E}$ defines a nondegenerate pairing on $\hat{A}$ with values in $B^*$, which is guaranteed by the metric double vector bundle structure; see \cite{Lean}. Coupled with the fact that the Courant bracket is closed on linear sections, one gets the following result.

\begin{pro}{\rm(\cite{Lean})}
The vector bundle $\hat{A}$ inherits a Courant algebroid structure with the pairing taking values in $B^*$, which  is called the \emph{fat} Courant algebroid of this $\VB$-Courant algebroid.
\end{pro}

Alternatively, we have
\begin{pro}\label{fat-E}
For a $\VB$-Courant algebroid $(\mathbb{E};A,B;M)$, its associated fat Courant algebroid is  a $B^*$-Courant algebroid.
\end{pro}
\begin{ex} {\rm{\bf (Standard $\VB$-Courant algebroid over a vector bundle)}\ \
For a vector bundle $E$, there is a standard $\VB$-Courant algebroid
\begin{equation*}
\xymatrix{ TE^*\oplus T^*E^*\ar@{->} [r]^{}\ar@ {->} [d]_{}
&E^* \ar@ {->} [d]^{}\\
{TM\oplus E} \ar@{->} [r]_{}&M }
\end{equation*}
with base $E^*$ and core $E^*\oplus T^*M \rightarrow M$. The corresponding exact sequence is given by
$$0\rightarrow \gl(E)\oplus T^*M\otimes E\longrightarrow \hat{A} \longrightarrow TM\oplus E\rightarrow 0.$$
Actually, by Example \ref{fat seq}, the corresponding fat Courant algebroid $\hat{A}$ here is exactly the omni-Lie algebroid  $\ol(E)=\dev E\oplus \jet E$. So the omni-Lie algebroid  is the linearization of the standard $\VB$-Courant algebroid. 
}
\end{ex}

\begin{ex} {\rm{\bf (Tangent $\VB$-Courant algebroid)}\ \ The tangent bundle $T\huaC$ of a Courant algebroid $\huaC\rightarrow M$
\begin{equation*}
\xymatrix{ T\huaC\ar@{->} [r]^{}\ar@ {->} [d]_{}
&TM \ar@ {->} [d]^{}\\
{\huaC} \ar@{->} [r]_{}&M }
\end{equation*}
carries a $\VB$-Courant algebroid structure with base $TM$ and core $\huaC\rightarrow M$.
The associated exact sequence is
$$ 0\rightarrow T^*M\otimes \huaC\longrightarrow \hat{\huaC} \longrightarrow \huaC\rightarrow 0.$$
Actually the fat Courant algebroid $\hat{\huaC}$ is $\jet \huaC$, which is a $T^*M$-Courant algebroid by Proposition \ref{fat-E}. So we get that on the jet bundle of a Courant algebroid, there is a $T^*M$-Courant algebroid structure. This result was first given in \cite{CLS2}.
}
\end{ex}

 A crossed module of Lie algebras consists of a pair of Lie algebras  $(\mathfrak{m},\mathfrak{g})$, an action $\triangleright$ of $\mathfrak{g}$ on $\mathfrak{m}$ and a Lie algebra morphism $\phi: \mathfrak{m}\rightarrow \mathfrak{g}$ such that
 $$ \phi(\xi)\triangleright \eta=[\xi,\eta]_\frkm,\quad \quad \phi(x\triangleright \xi)=[x,\phi(\xi)]_\g,\quad \forall ~x\in \mathfrak{g}, ~\xi, \eta\in \mathfrak{m}.$$

 Given a crossed module, there is an action $\rho: \mathfrak{g}\ltimes \mathfrak{g}^*\longrightarrow \frkX(\mathfrak{m}^*)$ of the natural quadratic Lie algebra $\mathfrak{g}\ltimes \mathfrak{g}^*$ on $\mathfrak{m}^*$
given by
 $$\rho(u+\alpha)=u\triangleright \cdot +\phi^*\alpha,$$
 where $u\triangleright \cdot\in \gl(\mathfrak{m})$  is viewed as a linear vector field on $\mathfrak{m}^*$ and $\phi^*\alpha\in \frkm^*$ is viewed as a constant vector field on $\mathfrak{m}^*$.
Note that  this action is coisotropic. We get the action Courant algebroid \cite{LM} $(\mathfrak{g}\ltimes \mathfrak{g}^*)\times \mathfrak{m}^*$ over $\mathfrak{m}^*$ with the anchor given by $\rho$ and the Dorfman bracket given by \begin{eqnarray}\label{action CA}
[e_1,e_2]=\huaL_{\rho(e_1)}e_2-\huaL_{\rho(e_2)}e_1+[e_1,e_2]_{\mathfrak{g}\ltimes \mathfrak{g}^*}+\rho^*\langle \dM e_1,e_2\rangle.
\end{eqnarray}
for any $e_1,e_2\in \Gamma((\mathfrak{g}\ltimes \mathfrak{g}^*)\times \mathfrak{m}^*)$. Here $\dM e_1\in \Omega^1(\mathfrak{m}^*,  \mathfrak{g}\ltimes \mathfrak{g}^*)$ is given by Lie derivatives: $(\dM e_1)(X)=\huaL_X e_1$ for $X\in \mathfrak{X}(\mathfrak{m}^*)$.
 Moreover, it is a $\VB$-Courant algebroid
\begin{equation*}
\xymatrix{ (\mathfrak{g}\ltimes \mathfrak{g}^*)\times \mathfrak{m}^*\ar@{->} [r]^{}\ar@ {->} [d]_{}
&\mathfrak{m}^* \ar@ {->} [d]^{}\\
\mathfrak{g} \ar@{->} [r]_{}&\ast }
\end{equation*}
with base $\mathfrak{m}^*$ and core $\mathfrak{g}^*$. See \cite{Li-Bland} for details. The associated exact sequence is
$$0\rightarrow \mathfrak{m}\otimes \mathfrak{g}^* \cong\Hom(\mathfrak{g},\mathfrak{m})\longrightarrow \hat{A}\longrightarrow \mathfrak{g}\rightarrow 0.$$
Since the double vector bundle is trivial,
we have $\hat{A}=\Hom(\mathfrak{g},\mathfrak{m})\oplus \mathfrak{g}$.

Moreover,  applying \eqref{action CA},   we get the Dorfman bracket on $
\Hom(\mathfrak{g},\mathfrak{m})\oplus \mathfrak{g}$.
\begin{pro}
With the above notations,
$(\Hom(\mathfrak{g},\mathfrak{m})\oplus \mathfrak{g},[\cdot,\cdot],(\cdot,\cdot)_\frkm,\rho=0 )$ is  an $\mathfrak{m}$-Courant algebroid, where the pairing $(\cdot,\cdot)_\frkm$ is given by $$( A+u,B+v)_\frkm=\half(Av+Bu),$$
and the Dorfman bracket is given by
\begin{eqnarray*}
[u,v]&=&[u,v]_\mathfrak{g};\\
{[A,B]}&=&A\circ \phi \circ B-B\circ \phi \circ A;\\
{[A,v]}&=&A\circ \ad_v^0-\ad_v^1\circ A+\cdot \triangleright Av+\phi(Av);\\
{[v,A]}&=&\ad_v^1\circ A-A\circ \ad_v^0
\end{eqnarray*}
for all $A,B\in \Hom(\mathfrak{g},\mathfrak{m}), u,v\in \mathfrak{g}.$ Here $\ad_v^0\in \gl(\mathfrak{g})$ and $\ad_v^1\in \gl(\mathfrak{m})$ are given by $\ad_v^0(u)=[v,u]_\mathfrak{g}$ and $\ad_v^1(a)=v\triangleright a$ respectively and $\cdot \triangleright Av\in \Hom(\mathfrak{g},\mathfrak{m})$ is defined by $(\cdot \triangleright Av)(u)=u\triangleright Av.$
\end{pro}
\pf By (\ref{action CA}), it is obvious that
\begin{eqnarray*}\label{equation1}
[u,v]=[u,v]_\mathfrak{g}.
\end{eqnarray*}
For $A,B\in \Hom(\mathfrak{g},\mathfrak{m}),v\in \mathfrak{g}$, applying (\ref{action CA}), we find
\begin{eqnarray*}\label{equation2}
[A,B]=\huaL_{\rho(A)} B-\huaL_{\rho(B)} A=\rho(A)B-\rho(B)A=A\circ \phi\circ B-B\circ \phi \circ A.
\end{eqnarray*}
Observe that $\huaL_{\rho(v)} A=\ad_v^1(A)=\ad_v^1\circ A$ and $[A,v]_{\mathfrak{g}\ltimes \mathfrak{g}^*}=-(\ad^0)_v^*A=A\circ \ad_v^0$. We have
\begin{eqnarray*}\label{equation3}
[A,v]\nonumber&=&\huaL_{\rho(A)} v-\huaL_{\rho(v)} A+[A,v]_{\mathfrak{g}\ltimes \mathfrak{g}^*}+\rho^*\langle \dM A,v\rangle\\ \nonumber&=&0-\ad_v^1\circ A+A\circ \ad_v^0+\cdot \triangleright Av+\phi(Av),
\end{eqnarray*}
where we have used
$$\rho^*\langle \dM A,v\rangle(u+B)=\rho(u+B)(Av)=u\triangleright Av+B(\phi(Av)).$$
Finally, we have
\begin{eqnarray*}\label{equation4}
[v,A]\nonumber&=&\huaL_{\rho(v)} A-\huaL_{\rho(A)} v+[v,A]_{\mathfrak{g}\ltimes \mathfrak{g}^*}+\rho^*\langle \dM v,A\rangle\\ &=&
\ad_v^1\circ A+0-A\circ \ad_v^0+0.
\end{eqnarray*}
This completes the proof.\qed

\begin{rmk}
This bracket can be viewed as a generalization of an omni-Lie algebra. See \cite[Example 5.2]{LSX} for more details.
\end{rmk}

 More generally, since the category of Lie 2-algebroids and the category of $\VB$-Courant algebroids are equivalent (see \cite{Li-Bland}), we get an $\E$-Courant algebroid from a Lie 2-algebroid. This construction first appeared in  \cite[Corollary 6.9]{Lean}.  Explicitly, let $(A_0\oplus A_{-1}, \rho_{A_0}, l_1, l_2=l_2^0+l_2^1,l_3)$ be a Lie 2-algebroid. Then we have an $A_{-1}$-Courant algebroid structure on \[ \Hom(A_0,A_{-1})\oplus A_0,\]
where the paring is given by  $$(D+u,D'+v)_{A_{-1}}=\half(Dv+D'u),\quad \quad~\forall~ D,D'\in \Gamma(\Hom(A_0,A_{-1})),~u,v\in \Gamma(A_0)$$ the anchor is
\[\rho: \Hom(A_0,A_{-1})\oplus A_0\to \dev A_{-1},\quad \quad \rho(D+u)=D\circ l_1+l_2^1(u,\cdot),\]
and the Dorfman bracket is given by
\begin{eqnarray*}
[u,v]&=&l_2^0(u,v)+l_3(u,v,\cdot);\\
{[D,D']}&=&D\circ l_1 \circ D'-D'\circ l_1 \circ D;\\
{[D,v]}&=&-l_2^1(v,D(\cdot))+D(l_2^0(v,\cdot))+l_2^1(\cdot, D(v))+l_1(D(v));\\
{[v,D]}&=&l_2^1(v,D(\cdot))-D(l_2^0(v,\cdot)).
\end{eqnarray*}

\section{Generalized complex structures on $\E$-Courant algebroids}
In this section, we introduce the notion of a generalized complex structure on an $\E$-Courant algebroid. We will see that it
unifies the usual generalized complex structure on an
even-dimensional manifold and the generalized contact structure
on an odd-dimensional manifold.
\begin{defi} \label{def-generalized-cplx}
A bundle map $\huaJ:\huaK\longrightarrow
\huaK$ is called a \emph{generalized almost complex structure} on an $\E$-Courant algebroid $(\huaK,\CDorfman{\cdot,\cdot},\pairE{\cdot,\cdot},\rho)$ if it satisfies the algebraic properties
 \begin{equation}\label{eqn:alg2}
 \huaJ^2=-1\quad\mbox{ and} \quad\pairE{\huaJ(X),\huaJ(Y)}=\pairE{X,Y}.
 \end{equation}
Furthermore, $\huaJ$ is called a  \emph{generalized  complex
structure} if  the following integrability condition is satisfied:
\begin{equation}\label{eqn:integral condition}
[\huaJ(X),\huaJ(Y)]_\huaK-[X,Y]_\huaK-\huaJ\big([\huaJ(X),Y]_\huaK+[X,\huaJ(Y)]_\huaK\big)=0,\quad\forall~X,~Y\in\Gamma(\huaK).
\end{equation}
\end{defi}

\begin{pro}\label{pro:J-J}
 Let $\huaJ:\huaK\longrightarrow
\huaK$ be a generalized almost complex structure on an $\E$-Courant algebroid $(\huaK,\CDorfman{\cdot,\cdot},\pairE{\cdot,\cdot},\rho)$. Then we
have $$\huaJ^\star|_{\huaK}=-\huaJ.$$
\end{pro}
\pf By \eqref{eqn:alg2}, for all $X,Y\in\Gamma(\huaK)$, we have
$$
\huaJ^\star(\huaJ(Y))(X)=\huaJ(Y)(\huaJ(X))=2\pairE{\huaJ(X),\huaJ(Y)}=2\pairE{X,Y}=Y(X).
$$
Since $X\in\Gamma(\huaK)$ is arbitrary, we have
$$
\huaJ^\star(\huaJ(Y))=Y,\quad\forall~Y\in\Gamma(\huaK).
$$
For any $Z\in\Gamma(\huaK)$, let $Y=-\huaJ (Z)$. By \eqref{eqn:alg2}, we
have $Z=\huaJ (Y)$.
 Then we have
$$
\huaJ^\star(Z)=\huaJ^\star(\huaJ (Y))=Y=-\huaJ (Z),
$$
which implies that $\huaJ^\star|_{\huaK}=-\huaJ.$ \qed

\begin{rmk}
  Generalized complex structures on  an $\E$-Courant algebroid $(\huaK,\CDorfman{\cdot,\cdot},\pairE{\cdot,\cdot},\rho)$ are in one-to-one correspondence with Dirac sub-bundles $S\subset \huaK\otimes\mathbb
 C$ such that $\huaK\otimes\mathbb
 C=S\oplus\overline{S}$. By a Dirac sub-bundle  of $\huaK$, we mean a sub-bundle
$S\subset \huaK$ which is   closed under the bracket
$\CDorfman{\cdot,\cdot}$ and satisfies $S=S^\bot$. The pair $(S, \overline{S})$ is  an $\E$-Lie bialgebroid in the sense of \cite{CLS2}.
 \end{rmk}

\begin{rmk}
Obviously, the notion of a generalized contact bundle associated to $L$, which was introduced in \cite{VitaglianoWade}, is a special case of Definition \ref{def-generalized-cplx} where $E$ is the line bundle $L$. In particular, if  $E$ is the trivial line bundle $L^\circ=M\times \mathbb R$,  we have
\begin{eqnarray*}
\dev L^\circ=TM\oplus\mathbb R,\quad \jet L^\circ=T^*M\oplus\mathbb R.
\end{eqnarray*}
Therefore, $\huaE^1(M)=\dev L^\circ\oplus \jet L^\circ$. Thus, a generalized complex structure on an $\E$-Courant algebroid unifies generalized complex structures on
even-dimensional manifolds and generalized contact bundles on odd-dimensional manifolds
\end{rmk}

\begin{ex}\label{ex:AV}{\rm
Consider the $\E$-Courant algebroid $A^*\otimes E \oplus A$
given in \cite[Example 2.9]{CLS2} for any Lie algebroid
$(A,[\cdot,\cdot]_A,a)$ and an $A$-module $E$. Twisted by a 3-cocycle $\Theta\in\Gamma(\wedge^3A^*,E)$, one obtains the AV-Courant algebroid introduced in \cite{AV courant} by Li-Bland.  Consider $\huaJ$ of the form
$\huaJ_D=\Big(\begin{array}{cc}-R_D&0\\0&D\end{array}\Big)$, where
$D\in\gl(A)$ and $R_D:A^*\otimes E\longrightarrow A^*\otimes E$ is
given by $
R_D(\phi)=\phi\circ D.
$
We get that $\huaJ$ is a generalized complex structure on the $\E$-Courant algebroid $A^*\otimes E \oplus A$ if and only if $D$ is a Nijenhuis operator on the Lie algebroid $A$ and $D^2=-1$.

Actually, $D^2=-1$ ensures that condition \eqref{eqn:alg2} holds. The Dorfman bracket on $\huaK=A^*\otimes E\oplus A$ is given by
\[[u+\Phi,v+\Psi]_\huaK=[u,v]_A+\huaL_u \Psi-\huaL_v \Phi+\rhowx \mathbbm{d}\Phi(v),\quad \forall u,v\in \Gamma(A),\Phi,\Psi\in \Gamma(A^*\otimes E),\]
where $\rhowx: \jet E\to A^*\otimes E$ is the dual of the $A$-action $\rho: A\to \dev E$ on $E$.
Then it is straightforward to see that the integrability condition \eqref{eqn:integral condition}  holds   if and only if $D$ is a Nijenhuis operator on $A$.
}

\end{ex}

Any generalized complex structure on a Courant algebroid induces a Poisson structure on the base manifold (see e.g. \cite{Barton}). Similarly,  any generalized complex structure on an $\E$-Courant algebroid induces a Lie algebroid or a local Lie algebra structure (\cite{KirillovLocal}) on $E$ .

\begin{thm}\label{thm:BaseLieAlgebroid}
Let $\huaJ:\huaK\longrightarrow
\huaK$ be  a generalized complex structure on an $\E$-Courant algebroid $(\huaK,\CDorfman{\cdot,\cdot},\pairE{\cdot,\cdot},\rho)$. Define a bracket operation $[\cdot,\cdot]_E:\Gamma(E)\wedge\Gamma(E)\longrightarrow\Gamma(E)$ by
\begin{equation}\label{eq:LiebracketJ}
 [u,v]_E\defbe 2\pairE{\huaJ \rho^\star\jetd u,\rho^\star\jetd v}=(\rho\circ\huaJ \circ\rho^\star)(\jetd u)(v),\quad \quad \forall ~u,~v\in \Gamma(E).
\end{equation}
Then $(E,[\cdot,\cdot]_E,\jd\circ \rho\circ \huaJ\circ \rho^\star\circ \jetd)$ is a Lie algebroid  when ${\rm rank}(E)\geq 2$ and $(E,[\cdot,\cdot]_E) $ is a local Lie algebra   when ${\rm rank}(E)=1$ .
 \end{thm}
  \pf
The bracket is obviously skew-symmetric.  By the integrability of $\huaJ$, we have
\[[\huaJ(\rho^\star\jetd u),\huaJ(\rho^\star\jetd v)]_\huaK-[\rho^\star\jetd u,\rho^\star\jetd v]_\huaK-\huaJ\big([\huaJ(\rho^\star\jetd u),\rho^\star\jetd v]_\huaK+[\rho^\star\jetd u,\huaJ(\rho^\star\jetd v)]_\huaK\big)=0.\]
Pairing with $\rho^\star\jetd w$ for $w\in \Gamma(E)$, by (\rm{EC-3}) in Definition \ref{defi: $E$-C A} and the first equation in Lemma \ref{[e,Df]}, we have
\begin{eqnarray*}
&&\pairE{[\huaJ(\rho^\star\jetd u),\huaJ(\rho^\star\jetd v)]_\huaK,\rho^\star\jetd w}\\ &=&\rho(\huaJ\rho^\star\jetd u)\pairE{\huaJ\rho^\star\jetd v,\rho^\star\jetd w}-\pairE{\huaJ\rho^\star\jetd v,[\huaJ\rho^\star\jetd u,\rho^\star\jetd w]_\huaK}\\ &=&
2\pairE{\rho^\star\jetd \pairE{\huaJ\rho^\star\jetd v,\rho^\star\jetd w},\huaJ\rho^\star\jetd u}-
2\pairE{\huaJ\rho^\star\jetd v,\rho^\star\jetd \pairE{\huaJ\rho^\star\jetd u,\rho^\star\jetd w}}\\ &=&\frac{1}{2}[u,[v,w]_E]_E-\frac{1}{2}[v,[u,w]_E]_E.
\end{eqnarray*}
 By (\rm{EC-1}) and (\rm{EC-5}) in Definition \ref{defi: $E$-C A}, we have
 $$
 \pairE{[\rho^\star\jetd u,\rho^\star\jetd v]_\huaK,\rho^\star\jetd w}=0.
 $$
Finally, using Lemma \ref{[e,Df]}, we have
\begin{eqnarray*}
&&\pairE{[\huaJ(\rho^\star\jetd u),\rho^\star\jetd v]_\huaK+[\rho^\star\jetd u,\huaJ(\rho^\star\jetd v)]_\huaK),
\huaJ \rho^\star\jetd w}\\ &=&2\pairE{\rho^\star\jetd \pairE{\huaJ \rho^\star\jetd u,\rho^\star\jetd v},\huaJ \rho^\star\jetd w}+0\\ &=&\frac{1}{2}[w,[u,v]_E]_E.
\end{eqnarray*}
Thus we get the Jacobi identity for $[\cdot,\cdot]_E$.
To see the Leibniz rule, by definition, we have
\[[u,fv]_E=f[u,v]_E+\jd \rho \huaJ\rho^\star\jetd (u)(f)v.\]
So it is a Lie algebroid structure if and only if $\jd \circ \rho \circ \huaJ\circ \rho^\star\circ \jetd: E\to TM$ is a bundle map, which is always true when ${\rm rank}(E)\geq 2$ (see the proof of \cite[Theorem 3.11]{CLomni}).
\qed

\section{Generalized complex structures on omni-Lie algebroids}

In  this section, we  study  generalized complex structures on the omni-Lie
algebroid  $\ol(E)$.
 We view $\ol(E)$   as a sub-bundle of $\Hom(\ol(E),E)$ by the
nondegenerate $E$-valued pairing $\pairE{\cdot,\cdot}$, i.e.
$$
e_2(e_1)\triangleq2\pairE{e_1,e_2},\quad \forall~ e_1, ~e_2\in\Gamma(\ol(E)).
$$

By Proposition \ref{pro:J-J}, we have

\begin{cor}\label{cor:equivalent}
 A bundle map $\huaJ:\ol(E)\longrightarrow \ol(E)$ is a generalized almost complex structure on the omni-Lie
algebroid $\ol(E)$ if and only if the following conditions are
satisfied:
 \begin{eqnarray}
  \huaJ^2=-\Id,\quad
 \huaJ^\star|_{\ol(E)}=-\huaJ.
 \end{eqnarray}
\end{cor}

Since $\ol(E)$ is the direct sum of $\dev E$ and $\jet E$, we can
write a generalized almost complex structure $\huaJ$ in the form of
a matrix. To do that, we need some preparations.

Vector bundles $\Hom(\wedge^k \dev E,E)_{\jet E}$ and $\Hom(\wedge^k
\jet E,E)_{\dev E}$ are introduced in \cite{CLS2} and
\cite{shengdeformation} to study deformations of omni-Lie algebroids
and deformations of Lie algebroids respectively. More precisely, we have
\begin{eqnarray*}
\Hom(\wedge^k \dev E,E)_{\jet E}&\triangleq &\set{\mu\in
\Hom(\wedge^k \dev E, E)~|~ \Img (\mu_\natural)\subset \jet E},
\quad (k\geq 2),\\
\Hom(\wedge^k\jet E,E)_{\dev E}& \triangleq &\set{\frkd\in
\Hom(\wedge^k \jet E, E)~|~ \Img (\frkd^\sharp)\subset \dev E},
\quad (k\geq 2),
\end{eqnarray*}
in which $\mu_\natural:\wedge^{k-1}\dev E\longrightarrow\Hom(\dev E,
E)$ is given by
\begin{equation*}
 \mu_\natural(\frkd_1,\cdots,
\frkd_{k-1})(\frkd_k)=\mu(\frkd_1,\cdots,
\frkd_{k-1},\frkd_k),\quad\forall ~\frkd_1,\cdots,\frkd_k\in\dev E,
\end{equation*}
and $\frkd^\sharp$ is defined similarly. By \eqref{eqn:jetE}, for
any $\mu\in\Hom(\wedge^k \dev E,E)_{\jet E}$, we have
\begin{equation}\label{eqn:jetkE}
\mu(\frkd_1,\cdots, \frkd_{k-1},\Phi)=\Phi\circ\mu(\frkd_1,\cdots,
\frkd_{k-1},\Id_E).
\end{equation}
Furthermore, $(\Gamma(\Hom(\wedge^\bullet\dev E,E)_{\jet E}),\jetd)$
is a subcomplex of
$(\Gamma(\Hom(\wedge^\bullet\dev E,E),\jetd)$, where $\jetd$ is the
coboundary operator of the gauge Lie algebroid $\dev E$
with the obvious action on   $E$.

\begin{pro}\label{pro:J}
Any generalized almost complex structure $\huaJ$ on the omni-Lie algebroid $\ol(E)$ must be of the form
\begin{equation}\label{eqn:type J}
\huaJ=\Big(\begin{array}{cc} N& \pi^\sharp\\
\sigma_\natural & -N^\star
\end{array}\Big),
\end{equation}
where $N:\dev E\longrightarrow\dev E$ is a bundle map satisfying
$N^\star(\jet E)\subset\jet E$, $\pi\in\Gamma(\Hom(\wedge^2\jet
E,E)_{\dev E})$, $~\sigma\in\Gamma(\Hom(\wedge^2\dev E,E)_{\jet E})$
such that the following conditions hold:
\begin{equation*}\pi^\sharp\circ\sigma_\natural+N^2=-\Id,\quad
N\circ\pi^\sharp=\pi^\sharp\circ N^\star,\quad \sigma_\natural\circ
    N=N^\star\circ\sigma_\natural.
    \end{equation*}
\end{pro}

\pf  By Corollary \ref{cor:equivalent}, for any generalized almost complex
structure $\huaJ$, we have $\huaJ^\star|_{\ol(E)}=-\huaJ$. Thus $\huaJ$
must be of the form
$$
\huaJ=\Big(\begin{array}{cc} N& \phi\\
\psi & -N^\star
\end{array}\Big),
$$
where  $N:\dev E\longrightarrow\dev E$  is a bundle map satisfying
$N^\star(\jet E)\subset\jet E$, $\phi:\jet E\longrightarrow\dev E$
and $\psi:\dev E\longrightarrow\jet E$ are bundle maps satisfying
$$
-\pairE{\phi(\mu),\nu}=\pairE{\mu,\phi(\nu)},\quad \quad -\pairE{\psi(\frkd),\frkt}=\pairE{\frkd,\psi(\frkt)}.
$$
Therefore,  we have $\phi=\pi^\sharp$, for some
$\pi\in\Gamma(\Hom(\wedge^2\jet E,E)_{\dev E})$, and
$\psi=\sigma_\natural$ for some $~\sigma\in\Gamma(\Hom(\wedge^2\dev
E,E)_{\jet E})$. This finishes the proof of the first part. As for
the second part, it is straightforward to see that the
conditions  follow from the fact that $\huaJ^2=-\Id$.\qed

\begin{rmk}
For a line bundle $L$, we have $\jet L=\Hom(\dev L,L)$ and $\dev
L=\Hom(\jet L,L)$. Therefore, the condition $N^\star(\jet
L)\subset\jet L$ always holds.
\end{rmk}

\begin{thm}\label{thm:Main thm}
A generalized almost complex structure $\huaJ$ given by
$(\ref{eqn:type J})$ is a generalized complex structure on the omni-Lie algebroid ${\ol(E)}$ if and only
if
\begin{itemize}
\item[\rm{(i)}] $\pi$ satisfies the
equation:
\begin{equation}\label{eqn:pi}
\pi^\sharp([\mu,\nu]_\pi)=[\pi^\sharp(\mu),\pi^\sharp(\nu)]_\dev,\quad\forall \mu,\nu\in\Gamma(\jet E),
\end{equation}
where the bracket $[\cdot,\cdot]_\pi$ on $\Gamma(\jet E)$ is defined by
\begin{equation}\label{eqn:pi bracket}
[\mu,\nu]_\pi\defbe \frkL_{\pi^\sharp(\mu)}\nu-\frkL_{\pi^\sharp(\nu)}\mu-
\jetd\ppE{\pi^\sharp(\mu),\nu}_E.
\end{equation}
\item[\rm{(ii)}] $\pi$ and $N$ are related by the following
formula:
\begin{eqnarray}\label{eqn:pi N}
N^\star([\mu,\nu]_\pi)&=&\frkL_{\pi^\sharp(\mu)}(N^\star(\nu))-\frkL_{\pi^\sharp(\nu)}(N^\star(\mu))-\jetd
\pi(N^\star(\mu),\nu).
\end{eqnarray}


\item[\rm{(iii)}]$N$ satisfies the condition:
\begin{eqnarray}\label{eqn:N}
T(N)(\frkd,\frkt)=\pi^\sharp (i_{\frkd\wedge\frkt}\jetd\sigma),\quad\forall \frkd,\frkt\in\Gamma(\dev E),
\end{eqnarray}
where $T(N)$ is the Nijenhuis tensor of $N$ defined by
$$
T(N)(\frkd,\frkt)=[N(\frkd),N(\frkt)]_\dev-N([N(\frkd),\frkt]_\dev+[\frkd,N(\frkt)]_\dev-N[\frkd,\frkt]_\dev).
$$

\item[\rm{(iv)}] $N$ and $\sigma$ are related by the following condition
\begin{equation}\label{eqn:Nsigma}
\jetd\sigma(N(\frkd),\frkt,\frkk)+\jetd\sigma(\frkd,N(\frkt),\frkk)+\jetd\sigma(\frkd,\frkt,N(\frkk))=\jetd\sigma_N(\frkd,\frkt,\frkk),\quad\forall \frkd,\frkt,\frkk\in\Gamma(\dev E),
\end{equation}
where $\sigma_N\in\Gamma(\Hom(\wedge^2\dev E, E)_{\jet E})$ is
defined by
$$
\sigma_N(\frkd,\frkt)=\sigma(N(\frkd),\frkt).
$$
\end{itemize}
\end{thm}
\emptycomment{$\phi$ and $i_N \phi$ are closed. In other words, the condition in quasi-Poisson Nijenhuis reduces to $di_N d\sigma=0$. Here because the cohomology is zero, so the closed form is exact. Is this connected?}
\emptycomment{Is it right that $\huaJ$ is a generalized complex structure if and only if $(\dev E)_N\oplus (\jet E)_{\pi}$ is an omni-Lie algebroid and $\huaJ: (\dev E)_N\oplus (\jet E)_{\pi}\longrightarrow \dev E \oplus \jet E$ is an Omni-Lie algebroid morphism.}
\pf Consider the integrability condition (\ref{eqn:integral
condition}). In fact, there are two equations since $\Gamma(\ol(E))$ has two components $\Gamma(\dev E)$ and $\Gamma(\jet E)$.  First let $e_1=\mu,~e_2=\nu$ be elements in
$\Gamma(\jet E)$, we have
$\huaJ(\mu)=\pi^\sharp(\mu)-N^\star(\mu)$,
$\huaJ(\nu)=\pi^\sharp(\nu)-N^\star(\nu)$ and
$\Dorfman{\mu,\nu}=0$. Therefore, we obtain
\begin{eqnarray*}
 && \Dorfman{\pi^\sharp(\mu)-N^\star(\mu),\pi^\sharp(\nu)-N^\star(\nu)}
  -\huaJ(\Dorfman{\pi^\sharp(\mu)-N^\star(\mu),\nu}+\Dorfman{\mu,\pi^\sharp(\nu)-N^\star(\nu)})\\
 & =&[\pi^\sharp(\mu),\pi^\sharp(\nu)]_\dev-\pi^\sharp(\Lied_{\pi^\sharp(\mu)}\nu-i_{\pi^\sharp(\nu)}\jetd \mu)\\
 &&+
 N^\star(\Lied_{\pi^\sharp(\mu)}\nu-i_{\pi^\sharp(\nu)}\jetd \mu)-\Lied_{\pi^\sharp(\mu)}N^\star(\nu)+i_{\pi^\sharp(\nu)}\jetd N^\star(\mu)\\
 &=&0.
\end{eqnarray*}
Thus we get conditions \eqref{eqn:pi} and \eqref{eqn:pi N}.

Then let $e_1=\frkd\in\Gamma(\dev E)$ and $e_2=\mu\in\Gamma(\jet E)$, we have
$\huaJ(e_1)=N(\frkd)+\sigma_\natural(\frkd)$ and
$\huaJ(e_2)=\pi^\sharp(\mu)-N^\star(\mu)$. Therefore, we
obtain
\begin{eqnarray*}
 &&
 \Dorfman{N(\frkd)+\sigma_\natural(\frkd),\pi^\sharp(\mu)-N^\star(\mu)}-\Dorfman{\frkd,\mu}
  -\huaJ(\Dorfman{N(\frkd)+\sigma_\natural(\frkd)),\mu}+\Dorfman{\frkd,\pi^\sharp(\mu)-N^\star(\mu)})\\
 ~& =& [N(\frkd),\pi^\sharp(\mu)]_\dev-N[\frkd,\pi^\sharp(\mu)]_\dev-\pi^\sharp(\Lied_{N(\frkd)}\mu-\Lied_\frkd
 N^\star(\mu))\\&&+
 N^\star(\Lied_{N(\frkd)}\mu-\Lied_\frkd N^\star(\mu))-\Lied_{N(\frkd)}N^\star(\mu)-i_{\pi^\sharp(\mu)}\jetd\sigma_\natural(\frkd)
-\Lied_\frkd\mu-\sigma_\natural[\frkd,\pi^\sharp(\mu)]_\dev \\
 &=&0.
\end{eqnarray*}
Thus we have
\begin{eqnarray}
\label{eqn:xxi1}
[N(\frkd),\pi^\sharp(\mu)]_\dev&=&N[\frkd,\pi^\sharp(\mu)]_\dev+\pi^\sharp(\Lied_{N(\frkd)}\mu-\Lied_\frkd
 N^\star(\mu)),\\
\label{eqn:xxi2}~ N^\star(\Lied_{N(\frkd)}\mu-\Lied_\frkd
 N^\star(\mu))&=&\Lied_{N(\frkd)}N^\star(\mu)+i_{\pi^\sharp(\mu)}\jetd\sigma_\natural(\frkd)+\Lied_\frkd\mu+\sigma_\natural[\frkd,\pi^\sharp(\mu)]_\dev.
\end{eqnarray}
We claim that \eqref{eqn:xxi1} is equivalent to \eqref{eqn:pi N}. In
fact,
applying \eqref{eqn:xxi1} to $\nu\in \Gamma(\jet E)$ and \eqref{eqn:pi N} to $\frkd\in \Gamma(\dev E)$, we get the same equality.
\emptycomment{\begin{eqnarray*}
-\ppE{\Lied_{N(\frkd)}\nu,\pi^\sharp(\mu)}_E =\ppE{N[\frkd,\pi^\sharp(\mu)]_\dev,\nu}_E+\ppE{\mu,[N(\frkd),\pi^\sharp(\nu)]_\dev}_E-\frkd\ppE{\mu,N\pi^\sharp(\nu)}_E+\ppE{\mu,N[\frkd,\pi^\sharp \nu]_\dev}_E.
\end{eqnarray*}
}

Next let $e_1=\frkd$ and $e_2=\frkt$ be elements in $\Gamma(\dev
E)$, we have $\huaJ(e_1)=N(\frkd)+\sigma_\natural(\frkd)$ and
$\huaJ(e_2)=N(\frkt)+\sigma_\natural(\frkt)$. Therefore, we have
\begin{eqnarray*}
 &&
 \Dorfman{N(\frkd)+\sigma_\natural(\frkd),N(\frkt)+\sigma_\natural(\frkt)}-[\frkd,\frkt]_\dev
  -\huaJ(\Dorfman{N(\frkd)+\sigma_\natural(\frkd),\frkt}+\Dorfman{\frkd,N(\frkt)+\sigma_\natural(\frkt)})\\
 ~& =& [N(\frkd),N(\frkt)]_\dev-[\frkd,\frkt]_\dev-N([N(\frkd),\frkt]_\dev+[\frkd,N(\frkt)]_\dev)-\pi^\sharp(\Lied_\frkd
 \sigma_\natural(\frkt)-i_\frkt\jetd\sigma_\natural(\frkd))\\
 &&+\Lied_{N(\frkd)}\sigma_\natural(\frkt)-i_{N(\frkt)}\jetd\sigma_\natural(\frkd)
-\sigma_\natural([N(\frkd),\frkt]_\dev+[\frkd,N(\frkt)]_\dev)+N^\star(\Lied_\frkd
 \sigma_\natural(\frkt)-i_\frkt\jetd\sigma_\natural(\frkd))\\
 &=&0.
\end{eqnarray*}
Thus we have
\begin{eqnarray}
\label{eqn:xy1}
[N(\frkd),N(\frkt)]_\dev-[\frkd,\frkt]_\dev-N([N(\frkd),\frkt]_\dev+[\frkd,N(\frkt)]_\dev)&=&\pi^\sharp(\Lied_\frkd
 \sigma_\natural(\frkt)-i_\frkt\jetd\sigma_\natural(\frkd)),\\
\label{eqn:xy2}
\sigma_\natural([N(\frkd),\frkt]_\dev+[\frkd,N(\frkt)]_\dev)-\Lied_{N(\frkd)}\sigma_\natural(\frkt)+i_{N(\frkt)}\jetd\sigma_\natural(\frkd)
&=&N^\star(\Lied_\frkd
 \sigma_\natural(\frkt)-i_\frkt\jetd\sigma_\natural(\frkd)).
\end{eqnarray}
We claim that \eqref{eqn:xxi2} and \eqref{eqn:xy1} are equivalent.
In fact, applying \eqref{eqn:xxi2} and \eqref{eqn:xy1} to $\frkt\in\Gamma(\dev E)$ and $\mu\in \Gamma(\jet E)$ respectively, we get the same equality
\begin{eqnarray*}
\ppE{[N(\frkd),N(\frkt)]_\dev-[\frkd,\frkt]_\dev-N([N(\frkd),\frkt]_\dev+[\frkd,N(\frkt)]_\dev),\mu}_E =
\frkd\ppE{\pi^\sharp \sigma_\sharp(\frkt),\mu}_E\\
+\ppE{\sigma_\sharp(\frkt),[\frkd,\pi^\sharp\mu]_\dev}_E+\frkt\ppE{\sigma_\sharp(\frkd),\pi^\sharp(\mu)}_E
-\pi^\sharp(\mu)\ppE{\sigma_\sharp(\frkd),\frkt}_E-\ppE{\sigma_{\sharp}(\frkd),[\frkt,\pi^\sharp(\mu)]_\dev}_E.
\end{eqnarray*}

By the equality $\pi^\sharp\circ\sigma_\natural+N^2=-\Id$ and
\eqref{eqn:xy1}, we have
$$
[N(\frkd),N(\frkt)]_\dev+N^2[\frkd,\frkt]_\dev-N([N(\frkd),\frkt]_\dev+[\frkd,N(\frkt)]_\dev)=\pi^\sharp(\Lied_\frkd
 \sigma_\natural(\frkt)-i_\frkt\jetd\sigma_\natural(\frkd)-\sigma_\natural[\frkd,\frkt]_\dev),
$$
which implies that $T(N)(\frkd,\frkt)=\pi^\sharp
(i_{\frkd\wedge\frkt}\jetd\sigma)$. Thus \eqref{eqn:xy1} is
equivalent to \eqref{eqn:N}.

At last, we consider condition \eqref{eqn:xy2}. Acting on an
arbitrary $\frkk\in\Gamma(\dev E)$, we have
\begin{eqnarray*}
  &&N(\frkd)\pe{\sigma_\natural(\frkt),\frkk}-\pe{\sigma_\natural(\frkt),[N(\frkd),\frkk]_\dev}+\pe{\sigma_\natural(\frkk),[N(\frkd),\frkt]_\dev}
-N(\frkt)\pe{\sigma_\natural(\frkd),\frkk}\\&&
  +\frkk\pe{\sigma_\natural(\frkd),\frkt}+\pe{\sigma_\natural(\frkd),[N(\frkt),\frkk]_\dev}
  +\pe{\sigma_\natural(\frkk),[\frkd,N(\frkt)]_\dev}+\frkd\pe{\sigma_\natural(\frkt),N(\frkk)}\\
  &&-\pe{\sigma_\natural(\frkt),[\frkd,N(\frkk)]_\dev}
  -\frkt\pe{\sigma_\natural(\frkd),N(\frkk)}+N(\frkk)\pe{\sigma_\natural(\frkd),\frkt}+\pe{\sigma_\natural(\frkd),[\frkt,N(\frkk)]_\dev}\\
  &=&\jetd
  \sigma(N(\frkd),\frkt,\frkk)+\frkt\sigma(N(\frkd),\frkk)-\frkk\sigma(N(\frkd),\frkt)+\sigma([\frkt,\frkk]_\dev,N(\frkd))\\
&&+\jetd
  \sigma(\frkd,N(\frkt),\frkk)-\frkd\sigma(N(\frkt),\frkk)-\sigma([\frkd,\frkk]_\dev,N(\frkt))\\
  &&+\jetd
  \sigma(\frkd,\frkt,N(\frkk))+\sigma([\frkd,\frkt]_\dev,N(\frkk))\\
  &=&0.
\end{eqnarray*}
Note that the following equality holds:
$$\sigma(\frkd,N(\frkt))=-\pe{\sigma_\natural(N(\frkt)),\frkd}=-\pe{N^\star(\sigma_\natural(\frkt)),\frkd}=-\pe{\sigma_\natural(\frkt),N(\frkd)}=\sigma(N(\frkd),\frkt).$$
Therefore, we have
$$
(i_N\jetd\sigma)(\frkd,\frkt,\frkk)=\jetd\sigma_N(\frkd,\frkt,\frkk),
$$
which implies that \eqref{eqn:xy2} is equivalent to \eqref{eqn:Nsigma}.
 \qed

\begin{rmk}
  Let $\huaJ=\Big(\begin{array}{cc} N& \pi^\sharp\\
\sigma_\natural & -N^\star
\end{array}\Big)$ be a generalized complex structure on the omni-Lie algebroid $\ol(E)$. Then $\pi$ satisfies  \eqref{eqn:pi}. On one hand, in \cite{CLomni}, the authors showed that such $\pi$ will give rise to a Lie bracket $[\cdot,\cdot]_E$ on $\Gamma(E)$ via
$$
[u,v]_E=\pi^\sharp(\jetd u)(v),\quad\forall~u,~v\in\Gamma(E).
$$
On the other hand, by Theorem \ref{thm:BaseLieAlgebroid}, the generalized   complex structure $\huaJ$  will also induce a Lie algebroid structure on $E$ by
$\eqref{eq:LiebracketJ}$. By the equality
$$
\pi^\sharp=\rho\circ \huaJ\circ \rhowx,
$$
these two Lie algebroid structures on $E$ are the same.
\end{rmk}

\begin{rmk}
 Recall that  any $b\in \Gamma(\Hom(\wedge^2\dev E, E)_{\jet E})$ defines a
transformation $e^b: \ol(E) \lon \ol(E)$ defined by
\begin{equation}\label{e b}
e^b\left(\begin{array}{c}\frkd\\\mu\end{array}\right)
=\left(\begin{array}{cc}\Id&0\\b_\natural&\Id\end{array}\right)\left(\begin{array}{c}\frkd\\\mu\end{array}\right)=\left(\begin{array}{c}\frkd\\\mu+i_\frkd
b\end{array}\right).
\end{equation}
Thus, $e^b$ is an automorphism of the omni-Lie algebroid $\ol(E)$ if and
only if $\jetd b=0$. In this case,  $e^b$ is called a
\emph{{$B$-field transformation}.}
Genuinely,  an automorphism of the omni-Lie algebroid $\ol(E)$ is just the
composition of an automorphism of the vector bundle $E$ and a
 $B$-field transformation. In fact $B$-field transformations map generalized complex structures on $\ol(E)$ into new  generalized complex structures as follows:
 \begin{equation}\label{e b transf}
\huaJ^b=\left(\begin{array}{cc}\Id&0\\b_\natural&\Id\end{array}\right) \circ \huaJ \circ \left(\begin{array}{cc}\Id&0\\ -b_\natural&\Id\end{array}\right) .
\end{equation}
\end{rmk}

\emptycomment{Actually,, from a general $\E$-Courant algebroid, we can get a Lie bracket on $\Gamma(E)$, which gives rise to $\pi\in \Hom(\wedge^2 \jet E, E)_{\dev E}$ satisfying $\pi^\sharp=\rho\circ \huaJ\circ \rho^\star$. In particular, in the omni-Lie algebroid case, this $\pi^\sharp$ is the element in the top right corner of $\huaJ$.}

\begin{ex}{\rm
Let $D: E\to E$ be a bundle map satisfying $D^2=-\Id$.  Define $R_D: \dev E\to \dev E$ by $R_D(\frkd)=\frkd\circ D$ and $\hat{D}: \jet E\to \jet E$ by $\hat{D}(\jetd  u)=\jetd (Du)$ for $u\in \Gamma(E)$. Then
\begin{equation*}
\huaJ=\Big(\begin{array}{cc} R_D& 0\\
0 & -\hat{D}
\end{array}\Big),
\end{equation*}
is a generalized complex structure on $\ol(E)$. In fact, since
\[\ppE{R_D^\star(\jetd u),\frkd}_E=\ppE{\jetd u,\frkd\circ D}_E=\frkd (D(u))=\ppE{\hat{D}(\jetd u), \frkd}_E,\]
we have $R_D^\star=\hat{D}$. It is direct to check that the Nijenhuis tensor $T(R_D)$ vanishes and the condition $D^2=-\Id$ ensures that $R_D^2=-\Id$.
}
\end{ex}

Let $\pi\in  \Gamma(\Hom(\wedge^2\jet E,E)_{\dev E})$ and suppose that the induced map
$\pi^\sharp: \jet E\to \dev E$ is an isomorphism of vector bundles. Then the rank of $E$ is $1$ or equals to the dimension of $M$. We denote by $(\pi^\sharp)^{-1}$ the inverse of $\pi^\sharp$  and by $\pi^{-1}$ the corresponding element in $\Gamma(\Hom(\wedge^2\dev E,E)_{\jet E})$.
\begin{lem}\label{lem:equivalent}
 With the above notations, the following two statements are equivalent:
  \begin{itemize}
    \item[\rm(i)] $\pi\in  \Gamma(\Hom(\wedge^2\jet E,E)_{\dev E})$ satisfies \eqref{eqn:pi};
    \item[\rm(ii)] $\pi^{-1}$ is closed, i.e. $\jetd \pi^{-1}=0$.
  \end{itemize}
\end{lem}
\pf The conclusion follows from the following equality:
$$
\ppE{\pi^\sharp([\mu,\nu]_\pi)-[\pi^\sharp(\mu),\pi^\sharp(\nu)]_\dev,\gamma}_E=-\jetd \pi^{-1}(\pi^\sharp(\mu),\pi^\sharp(\nu),\pi^\sharp(\gamma)),\quad \forall~\mu,~\nu,~\gamma\in\Gamma(\jet E),
$$
which can be obtained by straightforward computations.\qed\vspace{3mm}

Let $(E,[\cdot,\cdot]_E,a)$ be a Lie algebroid.  Define
$\pi^\sharp:\jet E\to \dev E$ by
\begin{equation}\label{eq:pi}
 \pi^\sharp(\jetd u)(\cdot)=[u,\cdot]_E, \quad \forall~u\in \Gamma(E).
 \end{equation}
Then $\pi^\sharp$ satisfies  \eqref{eqn:pi}. Furthermore, $(\jet E,[\cdot,\cdot]_\pi,\jd\circ \pi^\sharp)$ is a Lie algebroid, where the bracket $[\cdot,\cdot]_\pi$ is given by \eqref{eqn:pi bracket}.
 By Theorem \ref{thm:Main thm} and Lemma \ref{lem:equivalent}, we have
\begin{cor}
Let $(E,[\cdot,\cdot]_E,a)$ be a Lie algebroid such that the induced map $\pi^\sharp: \jet E\to \dev E$ is an isomorphism. Then
\begin{equation*}
\huaJ=\Big(\begin{array}{cc} 0& \pi^\sharp\\
-(\pi^\sharp)^{-1} & 0
\end{array}\Big),
\end{equation*}
is a generalized complex structure on $\ol(E)$.
\end{cor}

\begin{ex}\rm{
 Let $(TM,[\cdot,\cdot]_{TM},\Id)$ be the tangent Lie algebroid. Define $\pi^\sharp:\jet (TM)\to \dev (TM)$ by $$\pi^\sharp(\jetd u)=[u,\cdot]_{TM}.$$ Then $\pi^\sharp$ is an isomorphism. 
See  \cite[Corollary 3.9]{CLomni} for details.  Then
\begin{equation*}
\huaJ=\Big(\begin{array}{cc} 0& \pi^\sharp\\
-(\pi^\sharp)^{-1} & 0
\end{array}\Big),
\end{equation*}
is a generalized complex structure on the omni-Lie algebroid $\ol(TM)$. }
\end{ex}

\begin{ex}\rm{
Let $(M,\omega)$ be a symplectic manifold and  $(T^*M,[\cdot,\cdot]_{\omega^{-1}},(\omega^\sharp)^{-1})$  the associated natural Lie algebroid. Define
$\pi^\sharp: \jet (T^*M)\to \dev (T^*M)$ by $$ \pi^\sharp(\jetd u)=[u,\cdot]_{\omega^{-1}},$$ which is an isomorphism (see \cite[Corollary 3.10]{CLomni}).  Then
\begin{equation*}
\huaJ=\Big(\begin{array}{cc} 0& \pi^\sharp\\
-(\pi^\sharp)^{-1} & 0
\end{array}\Big),
\end{equation*}
is a generalized complex structure on the omni-Lie algebroid $\ol(T^*M)$.
}\end{ex}

 At the end of this section, we introduce the notion of an algebroid-Nijenhuis structure,  which can give rise to generalized complex structures on the omni-Lie algebroid $\ol(E)$.

\begin{defi}
Let $(E,[\cdot,\cdot]_E,a)$  be a Lie algebroid, $N:\dev
E\longrightarrow\dev E$  a Nijenhuis operator on the Lie algebroid $(\dev E,[\cdot,\cdot]_\dev,\jd)$ satisfying
$N^\star(\jet E)\subset\jet E$ and $\pi:\jet E\longrightarrow\dev E$
 given by \eqref{eq:pi}.  $N$ and $\pi$ are said to be
compatible if $$N\circ\pi^\sharp=\pi^\sharp\circ N^\star,\quad\mbox{and}\quad
C(\pi,N)=0,$$ where
\begin{equation}
                   C(\pi,N)(\mu,\nu)\triangleq [\mu,\nu]_{\pi_N}-([N^\star (\mu),\nu]_\pi+ [\mu,N^\star(\nu)]_\pi -
                   N^\star[\mu,\nu]_\pi),\quad \forall
                   ~\mu,~\nu\in\Gamma(\jet E).
                   \end{equation}
            Here $\pi_N\in\Gamma(\Hom(\wedge^2\jet E,E)_{\dev E})$  is given by
            $$
            \pi_N(\mu,\nu)=\ppE{\nu,N\pi^\sharp(\mu)}_E,\quad\forall~\mu,~\nu\in\Gamma(\jet E).
            $$    If $N$ and $\pi$ are compatible, we call the pair $(\pi,N)$ an
                   \emph{algebroid-Nijenhuis} structure on the Lie algebroid  $(E,[\cdot,\cdot]_E,a)$.
\end{defi}

The following lemma is straightforward, we omit the proof.
\begin{lem}\label{lem:PN}
Let  $(E,[\cdot,\cdot]_E,a)$ be a Lie algebroid, $\pi$ given by \eqref{eq:pi} and $N:\dev
E\longrightarrow\dev E$   a Nijenhuis structure. Then $(\pi,N)$ is an algebroid-Nijenhuis structure on the Lie algebroid  $(E,[\cdot,\cdot]_E,a)$  if and only if
$N\circ\pi^\sharp=\pi^\sharp\circ N^\star$ and
$$
N^\star[\mu,\nu]_\pi=\Lied_{\pi(\mu)}N^\star(\nu)-\Lied_{\pi(\nu)}N^\star(\mu)-\jetd
\pi(N^\star(\mu),\nu).
$$
 \end{lem}

By Theorem
\ref{thm:Main thm} and Lemma \ref{lem:PN}, we have

\begin{thm}\label{thm:LN and g C}
Let  $(E,[\cdot,\cdot]_E,a)$ be a Lie algebroid, $\pi$ given by \eqref{eq:pi} and $N:\dev
E\longrightarrow\dev E$   a Nijenhuis structure. Then the following statements are equivalent:
\begin{itemize}
\item[\rm(a)] $(\pi,N)$ is an algebroid-Nijenhuis   structure and
$N^2=-\Id;$
\item[\rm(b)] $\huaJ=\Big(\begin{array}{cc} N& \pi^\sharp\\
0 & -N^\star
\end{array}\Big)$ is a generalized complex structure on the
omni-Lie algebroid $\ol(E)$.
\end{itemize}
\end{thm}

\begin{rmk}
An interesting special case  is that where $E=L$ is a line bundle. Then $(\pi, N)$  becomes a \emph{Jacobi-Nijenhuis structure} on $M$. Jacobi-Nijenhuis structures are studied by Luca Vitagliano and the third author in \cite{VitaglianoWadeH}.  In this case, $\pi$ defines a Jacobi  bi-derivation  $\{\cdot,\cdot\}$ of $L$ (i.e. a skew-symmetric bracket which is a first order differential operator, hence a derivation, in each argument). Moreover, this bi-derivation is compatible with $N$ in the sense that   $\pi^\sharp \circ N^\star = N \circ \pi^\sharp$ and $C(\pi, N)=0$.  It defines a new Jacobi  bi-derivation $\{\cdot,\cdot\}_N$. Furthermore,
$(\{\cdot,\cdot\}, \{\cdot,\cdot\}_N)$ is a \emph{Jacobi bi-Hamiltonian} structure, i.e.~$\{\cdot,\cdot\}$, $\{\cdot,\cdot\}_N$ and $\{\cdot,\cdot\} + \{\cdot,\cdot\}_N$ are all Jacobi brackets.
\end{rmk}

\section{Generalized complex structures on omni-Lie algebras}
In this section, we consider the case that $E$ reduces to a vector
space $V$. Then we have
$$
\dev V=\gl(V),\quad \jet V=V.
$$
Furthermore, the pairing \eqref{Eqt:conpairingE} reduces
to
\begin{eqnarray}
\ppE{A,u}_V=Au, \quad \forall~A\in\gl(V),~u\in V.
\end{eqnarray}
Any $u\in V$ is a linear map from $\gl(V)$ to $V$,
$$
u(A)=\ppE{A,u}_V=Au.
$$
Therefore, an omni-Lie algebroid reduces to an omni-Lie algebra,
which is introduced by Weinstein in \cite{weinstein:omni} to study
the linearization of the standard Courant algebroid.

\begin{defi}An omni-Lie algebra associated to $V$ is a triple $(\gl(V)\oplus V,\Dorfman{\cdot,\cdot},\pairV{\cdot,\cdot} )$, where  $\pairV{\cdot,\cdot}$ is a nondegenerate
symmetric pairing  given by
 \begin{eqnarray}
\pairV{A+u,B+v}=\frac{1}{2}(Av+Bu), \quad \forall A,B\in\gl(V),~u,v\in V,
\end{eqnarray}
and $\Dorfman{\cdot,\cdot}$ is a
bracket operation given by
\begin{eqnarray}\label{eqn:bracketalgebra}
\Dorfman{A+u,B+v}&=&[A,B]+Av.
\end{eqnarray}
\end{defi}
We will simply denote an omni-Lie algebra associated to a vector space $V$ by $\ol(V)$.

\begin{lem}\label{lem:wedge2glV}
  For any vector space $V$, we have \begin{eqnarray*}
  \Hom(\wedge^2\gl(V),V)_V&=&0,\\\Hom(\wedge^2V,V)_{\gl(V)}&=&\Hom(\wedge^2V,V).
  \end{eqnarray*}
\end{lem}
\pf In fact, for any $\phi\in\Hom(\wedge^2\gl(V),V)_V$ and
$A, B\in\gl(V)$, by \eqref{eqn:jetkE}, we have
$$
\phi(A\wedge B)=B\circ\phi(A\wedge \Id_V)=-B\circ
A\circ\phi(\Id_V\wedge \Id_V)=0.
$$
Therefore, $\phi=0$, which implies that
$\Hom(\wedge^2\gl(V),V)_V=0$. The second equality is obvious. \qed

\begin{pro}
  Any generalized almost complex structure $\huaJ:\gl(V)\oplus V\longrightarrow\gl(V)\oplus
  V$ on the omni-Lie algebra $\ol(V) $ is of the following
  form
\begin{equation}\label{eqn:formJ-V}\Big(\begin{array}{cc}-R_D&\pi^\sharp\\0&D\end{array}\Big),\end{equation}
  where $\pi\in\Hom(\wedge^2V,V)$, $D\in\gl(V)$ satisfying
  $D^2=-\Id_V$ and $\pi(Du,v)=\pi(u,Dv)$, and
  $R_D:\gl(V)\longrightarrow\gl(V)$ is the right multiplication,
  i.e. $R_D(A)=A\circ  D$.
\end{pro}
\pf By Proposition \ref{pro:J} and Lemma \ref{lem:wedge2glV}, we can
assume that a generalized almost complex structure is of the form
$\Big(\begin{array}{cc}N&\pi^\sharp\\0&-N^\star\end{array}\Big),$
where $N:\gl(V)\longrightarrow\gl(V)$ satisfies $N^\star\in\gl(V)$
and $N^2=-\Id_{\gl(V)}$, and $\pi\in\Hom(\wedge^2V,V)$. Let
$D=-N^\star$, then we have
$$
(Dv)(A)=-N^\star(v)(A)=-v(N(A))=-N(A)v.
$$
On the other hand, we have $(Dv)(A)=ADv$, which implies that
$N(A)=-R_D(A)$. It is obvious that $N^2=-\Id_{\gl(V)}$ is equivalent
to that $D^2=-\Id_V$. The proof is finished. \qed

\begin{thm}
  A generalized almost complex structure $\huaJ:\gl(V)\oplus V\longrightarrow\gl(V)\oplus
  V$ on the omni-Lie algebra $\ol(V)$ given by
  \eqref{eqn:formJ-V} is a generalized complex structure if and only
  if\begin{itemize}
\item[\rm(i)] $\pi$ defines a Lie algebra structure $[\cdot,\cdot]_\pi$ on
$V$;

\item[\rm(ii)]  $D^2=-\Id_V$
 and $D[u,v]_\pi=[u,Dv]_\pi$ for $u,v\in V$.
  \end{itemize}
  Thus, a generalized complex structure on the omni-Lie algebra $\ol(V)$ gives rise to a complex Lie algebra structure on $V$.
\end{thm}

\pf By Theorem \ref{thm:Main thm}, we have
$$
[u,v]_\pi=\pi^\sharp(u)(v)=\pi(u,v).
$$
Condition \eqref{eqn:pi} implies that $[\cdot,\cdot]_\pi$ gives a Lie
algebra structure on $V$. Condition \eqref{eqn:pi N} implies that
$D[u,v]_\pi=[u,Dv]_\pi$. The other conditions are valid.

The conditions $D^2=-\Id_V$ and $D[u,v]_\pi=[u,Dv]_\pi$ say by definition  that $D$ is  a complex Lie algebra structure on  $(V,[\cdot,\cdot]_\pi). $ This finishes the
proof. \qed


\begin{thebibliography}{999}

\bibitem{Barton}
 J. Barton and M. Sti\'{e}non, Generalized complex submanifolds, \emph{Pacific J.  Math.} 236 (2008), no. 1, 23-44.

 \bibitem{conformalCA}
D. Baraglia,  Conformal Courant algebroids and orientifold T-duality, \emph{Int. J. Geom. Methods Mod. Phys.} 10 (2013), no. 2, 1250084, 35 pp.

 \bibitem{BCG}
H. Bursztyn, G. Cavalcanti and M. Gualtieri,  Reduction of Courant algebroids and generalized complex structures,  \emph{Adv. Math.} 211 (2007), no. 2, 726-765.








\bibitem{CLomni}
Z. Chen and Z. J. Liu,  Omni-Lie algebroids,  \emph{J. Geom. Phys.}  60 (2010), no. 5, 799-808.



\bibitem{CLS2}
Z. Chen and Z. J. Liu and Y. H. Sheng, $\E$-Courant algebroids, \emph{Int. Math. Res. Notices}  (2010), no. 22, 4334-4376.




\bibitem{CrainicGCS}
M. Crainic, Generalized complex structures and Lie brackets,
 \emph{Bull. Braz. Math. Soc. (N.S.)} 42 (2011), no. 4, 559-578.







\bibitem{Grabowski marmo2}
J. Grabowski  and G. Marmo,  The graded Jacobi algebras and (co)homology, \emph{J. Phys. A: Math. Gen.} 36 (2003), 161-81.

\bibitem{gualtieri}
M. Gualtieri,  Generalized complex geometry, \emph{Ann. of Math.} (2) 174 (2011), no. 1, 75-123.

\bibitem{hitchin}
N. J. Hitchin, Generalized Calabi-Yau manifolds, \emph{Q. J. Math.} 54 (2003), no. 3,
281-308.




\bibitem{Wade contact manifold and g}
D. Iglesias-Ponte and A. Wade, Contact manifold and generalized
complex structures, \emph{J. Geom. Phys.}  53  (2005),  249-258.

\bibitem{Lean}
M. Jotz Lean, $N$-manifolds of degree $2$ and metric double vector bundles,  arXiv:1504.00880.

\bibitem{KirillovLocal}
A. Kirillov, Local Lie algebras, \emph{{Russian Math. Surveys}}
31 (1976), 55-76.



\bibitem{LSX}
H. Lang, Y. Sheng and X. Xu, Nonabelian omni-Lie algebras, \emph{Banach Center Publications}  110 (2016), 167-176.




\bibitem{AV courant}
D. Li-Bland, $AV$-Courant algebroids and generalized CR
structures, \emph{Canad. J. Math.} 63 (2011), no. 4, 938-960.

\bibitem{Li-Bland}
D. Li-Bland, $\huaL\huaA$-Courant algebroids and their applications, thesis, University of Toronto, 2012, arXiv:1204.2796v1.
\bibitem{LM}
D. Li-Bland and E. Meinrenken, Courant algebroids and Poisson geometry, \emph{Int. Math. Res. Not. IMRN}, \textbf{11} (2009), 2106-2145.

\bibitem{LWXmani}%
Z. J. Liu, A. Weinstein and P. Xu, Manin triples for Lie
bialgebroids, \emph{J. Diff. Geom.}  45 (1997), 547-574.






\bibitem{Mkz:GTGA}
K. Mackenzie, \emph{General theories of Lie groupoids and Lie
algebroids}, Cambridge University Press, 2005.

\bibitem{ME}
K. Mackenzie, Ehresmann doubles and Drindel'd doubles for Lie algebroids and Lie bialgebroids,  \emph{J. Reine Angew. Math.} 658 (2011), 193-245.




\bibitem{Roytenbergphdthesis}
D. Roytenberg,
\newblock {Courant algebroids, derived brackets and even symplectic
  supermanifolds}, PhD thesis, UC Berkeley, 1999, arXiv:math.DG/9910078.




\bibitem{shengdeformation}
Y. Sheng, On deformations of Lie algebroids, \emph{Results. Math.} 62 (2012), 103-120.



\bibitem{SXquasi-N}
M. Sti\'{e}non and P. Xu, Poisson quasi-Nijenhuis manifolds,
\emph{Comm. Math. Phys.}  270 (2007), no. 3, 709-725.

\bibitem{SXreduction}
M. Sti\'{e}non and P. Xu, Reduction of generalized complex structures, \emph{J. Geom. Phys.} 58 (2008), no. 1, 105-121.


\bibitem{VitaglianoWade}
L. Vitagliano and A. Wade, Generalized contact bundles, \emph{C. R. Math. Acad. Sci. Paris}  {\bf 354} (2016), 313-317.

\bibitem{VitaglianoWadeH}
L. Vitagliano and A. Wade, Holomorphic Jacobi manifolds,  arXiv:1609.07737.




\bibitem{weinstein:omni}
A. Weinstein, Omni-Lie algebras, Microlocal analysis of the Schr${\rm\ddot{o}}$dinger equation and related topics (Japanese) (Kyoto, 1999), {\em S\=urikaisekikenky\=usho K\=oky\=uroku} 1176(2000), 95-102.


\end{thebibliography}
\end{document}